\documentclass[a4paper,oneside]{amsart}

\RequirePackage{amsmath}
\RequirePackage{bm}
\RequirePackage{amssymb}
\RequirePackage{upref}
\RequirePackage{amsthm}
\RequirePackage{enumerate}
\RequirePackage{pb-diagram}
\RequirePackage{amsfonts}
\RequirePackage[mathscr]{eucal}
\RequirePackage{verbatim}
\RequirePackage{xr}
\RequirePackage{graphicx}
\RequirePackage[pdftex]{floatflt}
\usepackage{calc}
\usepackage{xspace}
\usepackage{color}

\newcommand{\cf}{cf.\@\xspace}
\newcommand{\resp}{resp.\@\xspace}

\newcommand{\aev}{a.e.\@\xspace}

\newcommand{\al}{\alpha}
\newcommand{\bet}{\beta}
\newcommand{\ga}{\gamma}
\newcommand{\de}{\delta }
\newcommand{\e}{\epsilon}

\newcommand{\f}{\varphi}

\newcommand{\m}{\mu}
\newcommand{\n}{\nu}

\newcommand{\s}{\sigma}
\newcommand{\x}{\xi}

\newcommand{\C}{\varGamma}
\newcommand{\D}{\varDelta}
\newcommand{\F}{\varPhi}
\newcommand{\Lam}{\varLambda}

\newcommand{\di}[1]{#1\nobreakdash-\hspace{0pt}dimensional}

\newcommand{\fv}[2]{#1\hspace{0pt}_{|_{#2}}}

\newcommand{\so}{{\mc S_0}}

\newcommand{\const}{\tup{const}}
\newcommand{\ndash}{\nobreakdash--}

\newcommand{\msp[1]}[1]{\mspace{#1mu}}

\newcommand{\R}[1][n+1]{{\protect\mathbb R}^{#1}}

\newcommand{\N}{{\protect\mathbb N}}

\newcommand{\eR}{\stackrel{\lower1ex \hbox{\rule{6.5pt}{0.5pt}}}{\msp[3]\R[]}}
\newcommand{\eN}{\stackrel{\lower1ex \hbox{\rule{6.5pt}{0.5pt}}}{\msp[1]\N}}
\newcommand{\eO}{\stackrel{\lower1ex
\hbox{\rule{6pt}{0.5pt}}}{\msc O}}

\DeclareMathOperator{\graph}{graph}

\newcommand\ra{\rightarrow}

\newcommand\pa{\partial}

\newcommand{\un}{\infty}
\newcommand{\A}{\forall}

\newcommand{\set}[2]{\{\,#1\colon #2\,\}}
\newcommand{\uu}{\cup}
\newcommand{\ii}{\cap}
\newcommand{\uuu}{\bigcup}

\newcommand{\uud}{ \stackrel{\lower 1ex \hbox {.}}{\uu}}
\newcommand{\uuud}[1]{ \stackrel{\lower 1ex \hbox {.}}{\uuu_{#1}}}
\newcommand\su{\subset}

\newcommand\eS{\emptyset}
\newcommand{\sminus}[1][28]{\raise 0.#1ex\hbox{$\scriptstyle\setminus$}}

\newcommand{\abs}[1]{\lvert#1\rvert}

\newcommand{\norm}[1]{\lVert#1\rVert}

\newcommand{\spd}[2]{\protect\langle #1,#2\protect\rangle}

\newcommand\ch[3]{\varGamma_{#1#2}^#3}
\newcommand\cha[3]{{\bar\varGamma}_{#1#2}^#3}

\newcommand{\riem}[4]{R_{#1#2#3#4}}
\newcommand{\riema}[4]{{\bar R}_{#1#2#3#4}}

\newcommand{\tit}{\textit}

\newcommand{\tup}{\textup}

\newcommand{\mc}{\protect\mathcal}
\newcommand{\msc}{\protect\mathscr}

\providecommand{\bysame}{\makebox[3em]{\hrulefill}\thinspace}

\newcommand{\ci}{\cite}

\newcommand{\bt}{\begin{thm}}
\newcommand{\bl}{\begin{lem}}
\newcommand{\bc}{\begin{cor}}
\newcommand{\bd}{\begin{definition}}
\newcommand{\bpp}{\begin{prop}}
\newcommand{\br}{\begin{rem}}
\newcommand{\bn}{\begin{note}}
\newcommand{\be}{\begin{ex}}
\newcommand{\bes}{\begin{exs}}
\newcommand{\bb}{\begin{example}}
\newcommand{\bbs}{\begin{examples}}
\newcommand{\ba}{\begin{axiom}}

\newcommand{\et}{\end{thm}}
\newcommand{\el}{\end{lem}}
\newcommand{\ec}{\end{cor}}
\newcommand{\ed}{\end{definition}}
\newcommand{\epp}{\end{prop}}
\newcommand{\er}{\end{rem}}
\newcommand{\en}{\end{note}}
\newcommand{\ee}{\end{ex}}
\newcommand{\ees}{\end{exs}}
\newcommand{\eb}{\end{example}}
\newcommand{\ebs}{\end{examples}}
\newcommand{\ea}{\end{axiom}}

\newcommand{\bp}{\begin{proof}}
\newcommand{\ep}{\end{proof}}
\newcommand{\eps}{\renewcommand{\qed}{}\end{proof}}

\newcommand{\bal}{\begin{align}}

\newcommand{\bi}[1][1.]{\begin{enumerate}[\upshape #1]}
\newcommand{\bia}[1][(1)]{\begin{enumerate}[\upshape #1]}
\newcommand{\bin}[1][1]{\begin{enumerate}[\upshape\bfseries #1]}
\newcommand{\bir}[1][(i)]{\begin{enumerate}[\upshape #1]}
\newcommand{\bic}[1][(i)]{\begin{enumerate}[\upshape\hspace{2\cma}#1]}
\newcommand{\bis}[2][1.]{\begin{enumerate}[\upshape\hspace{#2\parindent}#1]}
\newcommand{\ei}{\end{enumerate}}

\newcommand\ndots{\raise 0.47ex \hbox {,}\hskip0.06em\cdots %
     \raise 0.47ex \hbox {,}\hskip0.06em} 

\newcommand{\q}{\quad}
\newcommand{\qq}{\qquad}

\newcommand\nd{\noindent}

\newskip\Csmallskipamount                                                
\Csmallskipamount=\smallskipamount
\newskip\Cmedskipamount
\Cmedskipamount=\medskipamount
\newskip\Cbigskipamount
\Cbigskipamount=\bigskipamount

\newcommand\cvs{\vspace\Csmallskipamount}   
\newcommand\cvm{\vspace\Cmedskipamount}

\newskip\csa
\csa=\smallskipamount

\newskip\cma
\cma=\medskipamount

\newskip\cba
\cba=\bigskipamount

\newdimen\spt
\newcommand\citem{\cvs\advance\itemno by
1{(\romannumeral\the\itemno})\hskip3pt}
\newcommand{\bitem}{\cvm\nd\advance\itemno by
1{\bf\the\itemno}\hspace{\cma}}

\newcount\itemno
\newcommand{\las}[1]{\label{S:#1}}

\newcommand{\lae}[1]{\label{E:#1}}
\newcommand{\lat}[1]{\label{T:#1}}
\newcommand{\lal}[1]{\label{L:#1}}
\newcommand{\lad}[1]{\label{D:#1}}
\newcommand{\lac}[1]{\label{C:#1}}

\newcommand{\rt}[1]{Theorem~\ref{T:#1}}
\newcommand{\rl}[1]{Lemma~\ref{L:#1}}
\newcommand{\rd}[1]{Definition~\ref{D:#1}}
\newcommand{\rc}[1]{Corollary~\ref{C:#1}}

\newcommand{\re}[1]{\eqref{E:#1}}

\newskip\thmskip
\thmskip=\parindent

\newskip\hsk
\newenvironment{hinw}{\labelsep=0pt\begin{list}{}{\labelsep=0pt\itemindent=0pt\labelwidth=0pt\leftmargin=\parindent\rightmargin=0pt\partopsep=\cba}%
\item\it\nopagebreak\nopagebreak}%
{\end{list}}

\newcommand\bh{\begin{hinw}}
\newcommand{\eh}{\end{hinw}}

\newtheoremstyle{normal}
  {\cba}
  {\cba}
  {}
  {\thmskip}
  {\bfseries}
  {.}
  {\hsk}
  {}

\newtheoremstyle{abschnitt}
  {\cba}
  {\cba}
  {}
  {\thmskip}
  {\bfseries}
  {.}
  {\hsk}
  {}

\newtheoremstyle{italic}
  {\cba}
  {\cba}
  {\itshape}
  {\thmskip}
  {\bfseries}
  {.}
  {\hsk}
  {}

\newtheoremstyle{aufgaben}
  {\cba}
  {\cba}
  {}
  {}
  {\normalsize\bfseries}
  {.}
  {\hsk}
  {}

\newtheoremstyle{break}
  {\cba}
  {\cba}
  {\itshape}
  {}
  {\bfseries}
  {.}
  {\newline}
  {}

\swapnumbers
\theoremstyle{italic}
\newtheorem{thm}[subsection]{Theorem}
\newtheorem{lem}[subsection]{Lemma}
\newtheorem{prop}[subsection]{Proposition}
\newtheorem{cor}[subsection]{Corollary}

\theoremstyle{normal}
\newtheorem{rem}[subsection]{Remark}
\newtheorem{definition}[subsection]{Definition}
\newtheorem{example}[subsection]{Example}
\newtheorem{examples}[subsection]{Examples}
\newtheorem{ex}[subsection]{Exercise}
\newtheorem{note}[subsection]{}
\newtheorem{axiom}[subsection]{Axiom}

\theoremstyle{aufgaben}
\newtheorem{exs}[subsection]{Exercises}

\swapnumbers

\numberwithin{equation}{section}
\numberwithin{figure}{section}

\newenvironment{textequation}[1][0.8]
{\begin{equation}
\begin{aligned}
\begin{minipage}{#1\linewidth}}
{\end{minipage}
\end{aligned}
\end{equation}
\ignorespacesafterend}

\newcommand{\btext}{\begin{textequation}}
\newcommand{\etext}{\end{textequation}}

\usepackage[german,english]{babel}
\usepackage{graphicx}
\RequirePackage{amsmath}
\RequirePackage{bm}
\RequirePackage{amssymb}
\RequirePackage{upref}
\RequirePackage{amsthm}
\RequirePackage{enumerate}
\RequirePackage{pb-diagram}
\RequirePackage{amsfonts}
\RequirePackage[mathscr]{eucal}
\RequirePackage{verbatim}
\RequirePackage{xr}
\RequirePackage{graphicx}
\RequirePackage[pdftex]{floatflt}
\usepackage{calc}
\usepackage{xspace}

\RequirePackage{upref}
\RequirePackage{amsthm}
\RequirePackage{enumerate}
\usepackage[mathscr]{eucal}

\usepackage{xr-hyper}

\usepackage{calc}

\newlength{\oddsidemarginlength}
\newlength{\topmarginlength}

\newcounter{numberoflines}
\newcounter{tempcc}
\oddsidemargin=\oddsidemarginlength
\evensidemargin=\oddsidemargin
\topmargin=\topmarginlength
\usepackage{hyperref}

\begin{document}

\flushbottom


\title{On the CMC foliation of future ends of a spacetime}

\author{Claus Gerhardt}
\address{Ruprecht-Karls-Universit\"at, Institut f\"ur Angewandte Mathematik,
Im Neuenheimer Feld 294, 69120 Heidelberg, Germany}
\email{gerhardt@math.uni-heidelberg.de}
\urladdr{http://www.math.uni-heidelberg.de/studinfo/gerhardt/}
\thanks{}

%
\date{\today}
%


\begin{abstract}
We consider spacetimes with compact Cauchy hypersurfaces and with Ricci tensor
bounded from below on the set of  timelike unit vectors, and prove that the results
known for spacetimes satisfying the timelike convergence condition, namely, foliation
by CMC hypersurfaces, are also valid in the present
situation, if corresponding further assumptions are satisfied. 

In addition we show that the volume of any sequence of spacelike hypersurfaces,
which run into the future singularity, decays to zero provided there exists a time
function covering a future end, such that the level hypersurfaces have non-negative
mean curvature and decaying volume.  
\end{abstract}

\maketitle

\tableofcontents

\setcounter{section}{0}
\section{Introduction}

Let $N$ be a $(n+1)$-dimensional spacetime with a compact Cauchy hypersurface,
so that $N$ is topologically a product, $N=I\times \so$, where $\so$ is a compact
Riemannian manifold and $I=(a,b)$ an interval. The metric in $N$ can then be
expressed in the form
\begin{equation}\lae{0.1}
d\bar s^2=e^{2\psi}(-(dx^0)^2+\s_{ij}(x^0,x)dx^idx^j);
\end{equation}
$x^0$ is the time function and $(x^i)$ are local coordinates for $\so$.

If $N$ satisfies a future mean curvature barrier condition and the timelike
convergence condition, then a future end $N_+=[a_0,b)$ can be foliated by constant
mean curvature (CMC) spacelike hypersurfaces and the mean curvature of the leaves
can be used as a new time function, \cf \cite{cg1,cg:foliation}. Moreover, one of
Hawking's singularity results implies that $N$ is future timelike incomplete with finite
Lorentzian diameter for the future end.

\cvm
In this paper we want to extend these results to the case when the Ricci tensor is
only bounded from below on the set of  timelike unit vectors
\begin{equation}\lae{0.2}
\bar R_{\al\bet}\nu^\al\nu^\bet\ge -\Lam\qq\A\,\spd\nu\nu=-1
\end{equation}
for some $\Lam\ge 0$, and in addition, we want to show that the volume of the CMC
leaves decays to zero, if the future singularity is approached.

\cvm
Our results can be summarized in

\bt\lat{0.1}
Suppose that in a future end $N_+$ of $N$ the Ricci tensor satisfies the estimate
\re{0.2}, and suppose that a future mean curvature barrier exists, \cf \rd{1.2}, then a
slightly smaller future end $\tilde N_+$ can be foliated by CMC spacelike
hypersurfaces, and there exists a smooth time function $x^0$ such that the slices
\begin{equation}
M_\tau=\{x^0=\tau\},\qq \tau_0<\tau<\un,
\end{equation}
have mean curvature $\tau$ for some $\tau_0>\sqrt{n\Lam}$. The precise value of
$\tau_0$ depends on the mean curvature of a lower barrier.
\et

\bt\lat{0.2}
Suppose that a future end $N_+=[a_0,b)$ of $N$ can be covered by a time function
$x^0$ such that the mean curvature of the slices $M_t=\{x^0=t\}$ is
non-negative and the volume of $M_t$ decays to zero
\begin{equation}
\lim_{t\ra b}\abs{M_t}=0,
\end{equation}
then the volume $\abs{M_k}$ of any sequence of spacelike achronal\,\footnote{A subset $M\su N$ is said to be achronal, if any timelike piecewise $C^1$-curve intersects $M$ at most once.} hypersurfaces $M_k$
that approach $b$, i.e.,
\begin{equation}
\lim_k\inf_{M_k}x^0=b,
\end{equation}
decays to zero. Thus, in case the additional conditions of \rt{0.1} are also satisfied,
the volume of the CMC hypersurfaces $M_\tau$ converges to zero
\begin{equation}
\lim_{\tau\ra\un}\abs{M_\tau}=0.
\end{equation}
\et

$N$ is also future timelike incomplete, if there is a compact spacelike hypersurface
$M$ with mean curvature $H$ satisfying
\begin{equation}
H\ge H_0>\sqrt{n\Lam},
\end{equation}
due to a result in \cite{galloway:cft}. 

\section{Notations and definitions}\las 1

The main objective of this section is to state the equations of Gau{\ss}, Codazzi,
and Weingarten for space-like hypersurfaces $M$ in a \di {(n+1)} Lorentzian
manifold
$N$.  Geometric quantities in $N$ will be denoted by
$(\bar g_{ \al \bet}),(\riema  \al \bet \ga \de)$, etc., and those in $M$ by $(g_{ij}), 
(\riem ijkl)$, etc. Greek indices range from $0$ to $n$ and Latin from $1$ to $n$; the
summation convention is always used. Generic coordinate systems in $N$ resp.
$M$ will be denoted by $(x^ \al)$ resp. $(\x^i)$. Covariant differentiation will
simply be indicated by indices, only in case of possible ambiguity they will be
preceded by a semicolon, i.e., for a function $u$ in $N$, $(u_ \al)$ will be the
gradient and
$(u_{ \al \bet})$ the Hessian, but e.g., the covariant derivative of the curvature
tensor will be abbreviated by $\riema  \al \bet \ga{ \de;\e}$. We also point out that
\begin{equation}
\riema  \al \bet \ga{ \de;i}=\riema  \al \bet \ga{ \de;\e}x_i^\e
\end{equation}
with obvious generalizations to other quantities.

Let $M$ be a \tit{spacelike} hypersurface, i.e., the induced metric is Riemannian,
with a differentiable normal $\n$ which is timelike.

In local coordinates, $(x^ \al)$ and $(\x^i)$, the geometric quantities of the
spacelike hypersurface $M$ are connected through the following equations
\begin{equation}\lae{1.2}
x_{ij}^ \al= h_{ij}\n^ \al
\end{equation}
the so-called \tit{Gau{\ss} formula}. Here, and also in the sequel, a covariant
derivative is always a \tit{full} tensor, i.e.,

\begin{equation}
x_{ij}^ \al=x_{,ij}^ \al-\ch ijk x_k^ \al+ \cha  \bet \ga \al x_i^ \bet x_j^ \ga.
\end{equation}
The comma indicates ordinary partial derivatives.

In this implicit definition the \tit{second fundamental form} $(h_{ij})$ is taken
with respect to $\n$.

The second equation is the \tit{Weingarten equation}
\begin{equation}
\n_i^ \al=h_i^k x_k^ \al,
\end{equation}
where we remember that $\n_i^ \al$ is a full tensor.

Finally, we have the \tit{Codazzi equation}
\begin{equation}
h_{ij;k}-h_{ik;j}=\riema \al \bet \ga \de\n^ \al x_i^ \bet x_j^ \ga x_k^ \de
\end{equation}
and the \tit{Gau{\ss} equation}
\begin{equation}
\riem ijkl=- \{h_{ik}h_{jl}-h_{il}h_{jk}\} + \riema  \al \bet\ga \de x_i^ \al x_j^ \bet
x_k^ \ga x_l^ \de.
\end{equation}

Now, let us assume that $N$ is a globally hyperbolic Lorentzian manifold with a
\tit{compact} Cauchy surface. 
$N$ is then a topological product $\R[]\times \mc S_0$, where $\mc S_0$ is a
compact Riemannian manifold, and there exists a Gaussian coordinate system
$(x^ \al)$, such that the metric in $N$ has the form 
\begin{equation}\lae{1.7}
d\bar s_N^2=e^{2\psi}\{-{dx^0}^2+\s_{ij}(x^0,x)dx^idx^j\},
\end{equation}
where $\s_{ij}$ is a Riemannian metric, $\psi$ a function on $N$, and $x$ an
abbreviation for the spacelike components $(x^i)$, see \ci{gr},
\ci[p.~212]{he:book}, \ci[p.~252]{grh}, and \ci[Section~6]{cg1}.
We also assume that
the coordinate system is \tit{future oriented}, i.e., the time coordinate $x^0$
increases on future directed curves. Hence, the \tit{contravariant} timelike
vector $(\x^ \al)=(1,0,\dotsc,0)$ is future directed as is its \tit{covariant} version
$(\x_ \al)=e^{2\psi}(-1,0,\dotsc,0)$.

Let $M=\graph \fv u\so$ be a spacelike hypersurface
\begin{equation}
M=\set{(x^0,x)}{x^0=u(x),\,x\in\mc S_0},
\end{equation}
then the induced metric has the form
\begin{equation}
g_{ij}=e^{2\psi}\{-u_iu_j+\s_{ij}\}
\end{equation}
where $\s_{ij}$ is evaluated at $(u,x)$, and its inverse $(g^{ij})=(g_{ij})^{-1}$ can
be expressed as
\begin{equation}\lae{1.10}
g^{ij}=e^{-2\psi}\{\s^{ij}+\frac{u^i}{v}\frac{u^j}{v}\},
\end{equation}
where $(\s^{ij})=(\s_{ij})^{-1}$ and
\begin{equation}\lae{1.11}
\begin{aligned}
u^i&=\s^{ij}u_j\\
v^2&=1-\s^{ij}u_iu_j\equiv 1-\abs{Du}^2.
\end{aligned}
\end{equation}
Hence, $\graph u$ is spacelike if and only if $\abs{Du}<1$.

The covariant form of a normal vector of a graph looks like
\begin{equation}
(\n_ \al)=\pm v^{-1}e^{\psi}(1, -u_i).
\end{equation}
and the contravariant version is
\begin{equation}
(\n^ \al)=\mp v^{-1}e^{-\psi}(1, u^i).
\end{equation}
Thus, we have
\br Let $M$ be spacelike graph in a future oriented coordinate system. Then, the
contravariant future directed normal vector has the form
\begin{equation}
(\n^ \al)=v^{-1}e^{-\psi}(1, u^i)
\end{equation}
and the past directed
\begin{equation}\lae{1.15}
(\n^ \al)=-v^{-1}e^{-\psi}(1, u^i).
\end{equation}
\er

In the Gau{\ss} formula \re{1.2} we are free to choose the future or past directed
normal, but we stipulate that we always use the past directed normal for reasons
that we have explained in \ci[Section 2]{cg:indiana}.

Look at the component $ \al=0$ in \re{1.2} and obtain in view of \re{1.15}

\begin{equation}\lae{1.16}
e^{-\psi}v^{-1}h_{ij}=-u_{ij}- \cha 000\mspace{1mu}u_iu_j- \cha 0j0
\mspace{1mu}u_i- \cha 0i0\mspace{1mu}u_j- \cha ij0.
\end{equation}
Here, the covariant derivatives are taken with respect to the induced metric of
$M$, and
\begin{equation}
-\cha ij0=e^{-\psi}\bar h_{ij},
\end{equation}
where $(\bar h_{ij})$ is the second fundamental form of the hypersurfaces
$\{x^0=\const\}$.

An easy calculation shows
\begin{equation}
\bar h_{ij}e^{-\psi}=-\tfrac{1}{2}\dot\s_{ij} -\dot\psi\s_{ij},
\end{equation}
where the dot indicates differentiation with respect to $x^0$.

\cvm
Finally, let us define what we mean by a future mean curvature barrier.
\bd\lad{1.2}
Let $N$ be a globally hyperbolic spacetime with compact Cauchy hypersurface $\so$
so that $N$ can be written as a topological product $N=\R[]\times \so$ and its
metric expressed as
\begin{equation}\lae{0.4}
d\bar s^2=e^{2\psi}(-(dx^0)^2+\s_{ij}(x^0,x)dx^idx^j).
\end{equation}
Here, $x^0$ is a globally defined future directed time function and $(x^i)$ are local
coordinates for $\so$.
$N$ is said to have a \tit{future mean curvature barrier} \resp \tit{past mean
curvature barrier}, if there are sequences $M_k^+$ \resp $M_k^-$ of closed
spacelike achronal hypersurfaces such that
\begin{equation}
\lim_{k\ra\un} \fv H{M_k^+}=\un \q\tup{\resp}\q \lim_{k\ra\un} \fv H{M_k^-}=-\un
\end{equation}
and
\begin{equation}
\limsup \inf_{M_k^+}x^0> x^0(p)\qq\A\,p\in N
\end{equation}
\resp
\begin{equation}
\liminf \sup_{M_k^-}x^0< x^0(p)\qq\A\,p\in N.
\end{equation}
\ed

\cvm
A future mean curvature barrier certainly represents a singularity, at least if $N$
satisfies \re{0.2}, because of the future timelike incompleteness, but these
singularities need not be crushing, \cf \cite[Introduction]{cg:imcf}.

\section{Proof of \rt{0.1}}

Let us start with some simple but very useful observations: If, for a given coordinate
system $(x^\al)$, the metric has the form \re{0.1}, then the coordinate slices
$M(t)=\{x^0=t\}$ can be looked at as a solution of the evolution problem
\begin{equation}\lae{2.1}
\dot x=-e^\psi \nu,
\end{equation}
where $\nu=(\nu^\al)$ is the past directed normal vector. The embedding
$x=x(t,\xi)$ is then given as $x=(t,x^i)$, where $(x^i)$ are local coordinates for
$\so$.

\cvm
From the equation \re{2.1} we can immediately derive evolution equations for the
geometric quantities $g_{ij}, h_{ij}, \nu$ and $H=g^{ij}h_{ij}$ of $M(t)$, \cf, e.g.,
\cite[Section~3]{cg:indiana}.

\cvm
To avoid confusion with notations  for the geometric
quantities of other hypersurfaces, we occasionally denote the induced metric and
second fundamental of  coordinate slices by $\bar g_{ij},\bar h_{ij}$ and $\bar H$.
Thus, the evolution equations
\begin{equation}\lae{2.2}
\dot{\bar g}_{ij}=-2 e^\psi \bar h_{ij}
\end{equation}
and
\begin{equation}\lae{2.3}
\dot{\bar H}=-\D e^\psi +(\abs{\bar A}^2+\bar R_{\al\bet}\nu^\al\nu^\bet) e^\psi
\end{equation}
are valid.

The last equation is closely related to the derivative of the mean curvature operator:
Let $M_0$ be a smooth spacelike hypersurface and consider in a tubular
neighbourhood $\mc U$ of $M_0$ hypersurfaces $M$ that can be written as graphs
over $M_0$, $M=\graph u$, in the corresponding normal Gaussian coordinate system.
Then the mean curvature of $M$ can be expressed as
\begin{equation}\lae{2.4}
H=-\D u+\bar H+v^{-2}u^iu^j\bar h_{ij},
\end{equation}
\cf equation \re{1.16}, and hence, choosing $u=\e\f$, $\f\in C^2(M_0)$, we deduce
\begin{equation}\lae{2.5}
\begin{aligned}
\frac d{d\e}\fv H {\e=0}&=-\D \f +\dot{\bar H}\f\\[\cma]
&= -\D\f +(\abs{\bar A}^2+\bar R_{\al\bet}\nu^\al\nu^\bet)\f.
\end{aligned}
\end{equation}

\cvm
Next we shall prove that CMC hypersurfaces are monotonically ordered, if the mean
curvatures are sufficiently large.

\bl\lal{2.1}
Let $M_i=\graph u_i$, $i=1,2$, be two spacelike hypersurfaces such that the \resp
mean curvatures $H_i$ satisfy
\begin{equation}\lae{2.6b}
H_1<H_2
\end{equation}
where $H_2$ is constant, $H_2=\tau_2$, and
\begin{equation}\lae{2.7b}
\sqrt{n\Lam}<\tau_2,
\end{equation}
 then there holds
\begin{equation}\lae{2.7}
u_1<u_2.
\end{equation}
\el

\bp
We first observe that the weaker conclusion
\begin{equation}\lae{2.10}
u_1\le u_2
\end{equation}
is as good as the strict inequality in \re{2.7}, in view of the maximum principle.

\cvm
Hence, suppose that \re{2.10} is not valid, so that
\begin{equation}\lae{2.11}
E(u_1)=\set{x\in\so}{u_2(x)<u_1(x)}\ne\eS.
\end{equation}

Then there exist points $p_i\in M_i$ such that
\begin{equation}
0<d_0=d(M_2,M_1)= d(p_2,p_1) =\sup \set{d(p,q)}{(p,q)\in M_2\times M_1},
\end{equation}
where $d$ is the Lorentzian distance function. Let $\f$ be a maximal geodesic from
$M_2$ to $M_1$  realizing this distance with endpoints $p_2$ and $p_1$, and
parametrized by arc length.

Denote by $\bar d$ the Lorentzian distance function to $M_2$, i.e., for $p\in
I^+(M_2)$
\begin{equation}
\bar d(p)=\sup_{q\in M_2}d(q,p).
\end{equation}

Since $\f$ is maximal, $\C=\set{\f(t)}{0\le t<d_0}$ contains no focal points of
$M_2$,
\cf \cite[Theorem 34, p. 285]{bn}, hence there exists an open neighbourhood $\mc
V=\mc V(\C)$ such that $\bar d$ is smooth in $\mc V$, \cf \cite[Proposition
30]{bn}, because $\bar d$ is a component of the inverse of the normal exponential
map of $M_2$.

\cvm
Now, $M_2$ is the level set $\{\bar d=0\}$, and the level sets
\begin{equation}
M(t)=\set{p\in \mc V}{\bar d(p)=t}
\end{equation}
are  smooth hypersurfaces; $x^0=\bar d$ is a time function in $\mc V$ and
generates a normal Gaussian coordinate system, since $\spd{D\bar d}{D\bar d}=-1$.
Thus, the mean curvature $\bar H(t)$ of $M(t)$ satisfies the equation
\begin{equation}
\dot {\bar H}=\abs{\bar A}^2+\bar R_{\al\bet}\nu^\al\nu^\bet,
\end{equation}
\cf \re{2.3}, and therefore we have
\begin{equation}\lae{2.15b}
\dot{\bar H}\ge \tfrac1n\abs{\bar H}^2-\Lam >0,
\end{equation}
in view of \re{2.7b}.

\cvm
Next, consider a tubular neighbourhood $\mc U$ of $M_1$ with corresponding
normal Gaussian coordinates $(x^\al)$. The level sets
\begin{equation}
\tilde M(s)=\{x^0=s\},\qq-\e<s<0,
\end{equation}
lie in the past of $M_1=\tilde M(0)$ and are smooth for small $\e$.

Since the geodesic $\f$ is normal to $M_1$, it is also normal to $\tilde M(s)$ and
the length of the geodesic segment of $\f$ from $\tilde M(s)$ to $M_1$ is exactly
$-s$, i.e., equal to the distance from $\tilde M(s)$ to $M_1$, hence we deduce
\begin{equation}
d(M_2,\tilde M(s))=d_0+s,
\end{equation}
i.e., $\set{\f(t)}{0\le t\le d_0+s}$ is also a maximal geodesic from $M_2$ to $\tilde
M(s)$, and we conclude further that, for fixed $s$, the hypersurface $\tilde
M(s)\ii\mc V$ is contained in the past of $M(d_0+s)$ and touches $M(d_0+s)$ in
$p_s=\f(d_0+s)$. The maximum principle then implies
\begin{equation}
\fv H{\tilde M(s)}(p_s)\ge \fv H {M(d_0+s)}(p_s)>\tau_2,
\end{equation}
in view of \re{2.15b}.

On the other hand, the mean curvature of $\tilde M(s)$ converges to the mean
curvature of $M_1$, if $s$ tends to zero, hence we conclude
\begin{equation}
H_1(\f(d_0))\ge \tau_2,
\end{equation}
contradicting \re{2.6b}.
\ep

As an immediate conclusion we obtain

\bc\lac{2.2}
The CMC hypersurfaces with mean curvature
\begin{equation}
\tau>\sqrt{n\Lam}
\end{equation}
are uniquely determined.
\ec

\bp
Let $M_i=\graph u_i$, $i=1,2$, be two hypersurfaces with mean curvature $\tau$
and suppose, e.g., that
\begin{equation}
\set{x\in\so}{u_1(x)<u_2(x)}\ne\eS.
\end{equation}
Consider a tubular neighbourhood of $M_1$ with a corresponding future oriented
normal Gaussian coordinate system $(x^\al)$. Then the evolution of the mean
curvature of the coordinate slices satisfies
\begin{equation}
\dot{\bar H}=\abs{\bar A}^2+\bar R_{\al\bet}\nu^\al\nu^\bet \ge \frac
1n\abs{\bar H}^2-\Lam >0
\end{equation}
in a neighbourhood of $M_1$, i.e., the coordinate slices $M(t)=\{x^0=t\}$, with
$t>0$, have all mean curvature $\bar H(t)>\tau$. Using now $M_1$ and $M(t)$,
$t>0$, as barriers, we infer that for any $\tau'\in\R[]$, $\tau<\tau' <\bar H(t)$,
there exists a spacelike hypersurface $M_{\tau'}$ with mean curvature $\tau'$, such
that
$M_{\tau'}$ can be expressed as a graph over $M_1$, $M_{\tau'}=\graph u$, where
\begin{equation}
0<u<t.
\end{equation}

For a proof see \ci[Section~6]{cg1}; a different more transparent proof of this result
has been given in \cite{cg:mz}.

Writing $M_{\tau'}$ as graph over $\so$ in the original coordinate system without
changing the notation for $u$, we obtain
\begin{equation}\lae{2.14}
u_1<u,
\end{equation}
and, by choosing $t$ small enough, we may also conclude that
\begin{equation}\lae{2.15}
E(u)=\set{x\in\so}{u(x)<u_2(x)}\ne \eS,
\end{equation}
which is impossible, in view of the preceding result.
\ep

\bl
Under the assumptions of \rt{0.1}, let $M_{\tau_0}=\graph u_{\tau_0}$ be a CMC
hypersurface with mean curvature $\tau_0>\sqrt{n\Lam}$, then the future of
$M_{\tau_0}$ can be foliated by CMC hypersurfaces
\begin{equation}\lae{2.28}
I^+(M_{\tau_0})=\uuu_{\tau_0<\tau<\un}M_\tau.
\end{equation}
The $M_\tau$ can be written as graphs over $\so$
\begin{equation}
M_\tau=\graph u(\tau,\cdot),
\end{equation}
such that $u$ is strictly monotone increasing with respect to $\tau$, and continuous
in $[\tau_0,\un)\times \so$.
\el

\bp
The monotonicity and continuity of $u$ follows from \rl{2.1} and \rc{2.2}, in view of
the a priori estimates.

\cvm
Thus, it remains to verify the relation \re{2.28}. Let $p=(t,y^i)\in I^+(M_{\tau_0})$,
then we have to show $p\in M_\tau$ for some $\tau>\tau_0$.

In \cite[Theorem 6.3]{cg1} it is proved that there exists a family of CMC
hypersurfaces $M_\tau$
\begin{equation}
\set{M_\tau}{\tau_0\le \tau<\un},
\end{equation}
if there is a future mean curvature barrier.

Define $u(\tau,\cdot)$ by
\begin{equation}
M_\tau=\graph u(\tau,\cdot),
\end{equation}
then we have
\begin{equation}
u(\tau_0,y)<t<u(\tau^*,y)
\end{equation}
for some large $\tau^*$, because of the mean curvature barrier condition, which,
together with \rl{2.1}, implies that the CMC hypersurfaces run into the future
singularity, if $\tau$ goes to infinity.

In view of the continuity of $u(\cdot,y)$ we conclude that there exists $\tau_1>\tau_0$  such that
\begin{equation}
u(\tau_1,y)= t,
\end{equation}
 hence
$p\in M_{\tau_1}$.
\ep

\br
The continuity and monotonicity of $u$ holds in any coordinate system $(x^\al)$,
even in those that do not cover the future completely like the normal Gaussian
coordinates associated with a spacelike hypersurface, which are defined in a tubular
neighbourhood.
\er

The proof of \rt{0.1} is now almost finished. The remaining arguments are identical to
those in \cite[Section 2]{cg:foliation}, but for the convenience of the reader, we shall
briefly summarize the main steps.

\cvm
 We have to show that the mean curvature parameter $\tau$ can be
used as a time function in $\{\tau_0<\tau<\un\}$, i.e., $\tau$ should be smooth with
a non-vanishing gradient. Both properties are local properties.

\cvm
\tit{First step}: Fix an arbitrary $\tau'\in (\tau_0,\un)$, and consider a tubular
neighbourhood $\mc U$ of $M'=M_{\tau'}$. The $M_\tau\su \mc U$ can then be
written as graphs over $M'$, $M_\tau=\graph u(\tau,\cdot)$. For small $\e>0$ we
have
\begin{equation}
M_\tau\su \mc U\qq\A\,\tau\in (\tau'-\e,\tau'+\e)
\end{equation}
and with the help of the implicit function theorem we shall show that $u$ is smooth.
Indeed, define the operator $G$
\begin{equation}
G(\tau,\f)=H(\f)-\tau,
\end{equation}
where $H(\f)$ is an abbreviation for the mean curvature of $\graph \fv \f{M'}$. Then
$G$ is smooth and from \re{2.5} we deduce that $D_2G(\tau',0)\f$ equals
\begin{equation}
-\D\f+(\norm{ A}^2+\bar R_{\al\bet}\nu^\al\nu^\bet)\f,
\end{equation}
where the Laplacian, the second fundamental form and the normal correspond to $M'$.
Hence $D_2G(\tau',0)$ is an isomorphism and the implicit function theorem implies
that $u$ is smooth.

\cvm
\tit{Second step}: Still in the tubular neighbourhood of $M'$, define the coordinate
transformation
\begin{equation}
\F(\tau,x^i)=(u(\tau,x^i),x^i);
\end{equation}
note that $x^0=u(\tau,x^i)$. Then we have
\begin{equation}
\det D\F=\frac{\pa u}{\pa \tau}=\dot u.
\end{equation}

$\dot u$ is non-negative; if it were strictly positive, then $\F$ would be a
diffeomorphism, and hence $\tau$ would be smooth with non-vanishing gradient. A
proof, that $\dot u>0$, is given in \cite[Lemma 2.2]{cg:foliation}, but let us give a
simpler proof: The CMC hypersurfaces in $\mc U$ satisfy an equation
\begin{equation}
H(u)=\tau,
\end{equation}
where the left hand-side can be expressed as in \re{2.4}. Differentiating both sides
with respect to $\tau$ and evaluating for $\tau=\tau'$, i.e., on $M'$, where
$u(\tau',\cdot)=0$, we get
\begin{equation}
-\D\dot u+(\abs{A}^2+\bar R_{\al\bet}\nu^\al\nu^\bet)\dot u=1.
\end{equation}

In a point, where $\dot u$ attains its minimum, the maximum principle implies
\begin{equation}
(\abs{A}^2+\bar R_{\al\bet}\nu^\al\nu^\bet)\dot u\ge1,
\end{equation}
hence $\dot u\ne0$ and $\dot u$ is therefore strictly positive.

\section{Proof of \rt{0.2}}

Let $x^0$ be time function satisfying the assumptions of \rt{0.2}, i.e.,
$N_+=\{a_0<x^0<b\}$, the mean curvature of the slices $M(t)=\{x^0=t\}$ is
non-negative, and
\begin{equation}
\lim_{t\ra b}\abs{M(t)}=0,
\end{equation}
and let $M_k$ be a sequence of spacelike hypersurfaces such that
\begin{equation}
\lim\inf_{M_k}x^0=b.
\end{equation}

Let us write $M_k=\graph u_k$ as graphs over $\so$. Then
\begin{equation}
g_{ij}=e^{2\psi}(u_iu_j+\s_{ij}(u,x))
\end{equation}
is the induced metric, where we dropped the index $k$ for better readability, and the
volume element of $M_k$ has the form
\begin{equation}
d\m=v\sqrt{\det(\bar g_{ij}(u,x))}\,dx,
\end{equation}
where
\begin{equation}\lae{3.5}
v^2=1-\s^{ij}u_iu_j<1,
\end{equation}
and $(\bar g_{ij}(t,\cdot))$ is the metric of the slices $M(t)$.

From \re{2.2} we deduce
\begin{equation}\lae{3.6}
\frac d{dt}\sqrt{\det(\bar g_{ij}(t,\cdot))}=-e^\psi\bar H\sqrt{\det(\bar g_{ij})}\le
0.
\end{equation}

\cvm
Now, let $a_0<t<b$ be fixed, then for \aev $k$ we have
\begin{equation}\lae{3.7}
t<u_k
\end{equation}
and hence
\begin{equation}
\begin{aligned}
\abs{M_k}&=\int_\so v \sqrt{\det(\bar g_{ij}(u_k,x))}\,dx\\[\cma]
&\le \int_\so\sqrt{\det(\bar g_{ij}(t,x)}\,dx =\abs{M(t)},
\end{aligned}
\end{equation}
in view of \re{3.5}, \re{3.6} and \re{3.7}, and we conclude
\begin{equation}
\limsup\abs{M_k}\le \abs{M(t)}\qq\A\,a_0<t<b,
\end{equation}
and thus
\begin{equation}
\lim\abs{M_k}=0.
\end{equation}

\nocite{eh1}

\bibliographystyle{amsplain}

\begin{thebibliography}{10}

\bibitem{galloway:cft}
Lars Andersson and Gregory Galloway, \emph{Ds/cft and spacetime topology}, Adv.
  Theor. Math. Phys. \textbf{68} (2003), 307--327, 17 pages,
  \href{http://arXiv.org/pdf/hep-th/0202161}{hep-th/0202161}.
  
\bibitem{eh1}
Klaus Ecker and Gerhard Huisken, \emph{{Parabolic methods for the construction
  of spacelike slices of prescribed mean curvature in cosmological
  spacetimes.}}, Commun. Math. Phys. \textbf{135} (1991), no.~3, 595--613.

\bibitem{cg1}
Claus Gerhardt, \emph{H-surfaces in {{L}orentzian} manifolds}, Commun. Math.
  Phys. \textbf{89} (1983), 523\ndash 553.

\bibitem{cg:indiana}
\bysame, \emph{Hypersurfaces of prescribed curvature in {L}orentzian
  manifolds}, Indiana Univ. Math. J. \textbf{49} (2000), 1125\ndash 1153,
  \href{http://www.math.uni-heidelberg.de/studinfo/gerhardt/bibtexcg00b.html}{%
pdf file}.

\bibitem{cg:mz}
\bysame, \emph{Hypersurfaces of prescribed mean curvature in {L}orentzian
  manifolds}, Math. Z. \textbf{235} (2000), 83\ndash 97.

\bibitem{cg:foliation}
\bysame, \emph{On the foliation of space-time by constant mean curvature
  hypersurfaces}, 2003, e-print, 7 pages,
  \href{http://arXiv.org/pdf/math.DG/0304423}{math.DG/0304423}.

\bibitem{cg:imcf}
\bysame, \emph{The inverse mean curvature flow in cosmological spacetimes},
  2004, 24 pages, \href{http://arxiv.org/pdf/math.DG/0403097}{math.DG/0403097}.

\bibitem{gr}
R.P. Geroch, \emph{The domain of dependence}, J. Math. Phys. \textbf{11}
  (1970), 437\ndash 449.

\bibitem{grh}
R.P. Geroch and G.T. Horowitz, \emph{Global structure of spacetime}, General
  Relativity. An Einstein centenary survey (S.W. Hawking and W.~Israel, eds.),
  Cambridge University Press, 1973, p.~212\ndash293.

\bibitem{he:book}
S.~W. Hawking and G.~F.~R. Ellis, \emph{The large scale structure of
  space-time}, Cambridge University Press, London, 1973. \MR{54 \#12154}

\bibitem{bn}
Barrett O'Neill, \emph{{Semi-Riemannian geometry. With applications to
  relativity.}}, {Pure and Applied Mathematics, 103. New York-London etc.:
  Academic Press. XIII}, 1983.

\end{thebibliography}
\providecommand{\bysame}{\leavevmode\hbox to3em{\hrulefill}\thinspace}
\providecommand{\MR}{\relax\ifhmode\unskip\space\fi MR }
\providecommand{\MRhref}[2]{%
  \href{http://www.ams.org/mathscinet-getitem?mr=#1}{#2}
}
\providecommand{\href}[2]{#2}



\end{document}